\documentclass[12pt]{article}
\usepackage[]{amsmath,amssymb,amsfonts,latexsym,amsthm,enumerate}
\newtheorem{thm}{Theorem}
\newtheorem*{thm*}{Theorem}
\newtheorem{lemma}[thm]{Lemma}
\newtheorem*{lemma*}{Lemma}

\newtheorem*{prop*}{Proposition}
\newtheorem{cor}[thm]{Corollary}
\theoremstyle{remark}
\newtheorem*{rmk*}{Remark}
\newtheorem*{rmks*}{Remarks}
\newtheorem*{not*}{Notation}
\newtheorem*{claim*}{Claim}
\newtheorem*{fact*}{Fact}

\theoremstyle{definition}
\newtheorem{dfn}{Definition}
\newtheorem*{namedprop}{\theoremname}
\newcommand{\theoremname}{testing}
\newenvironment{lemma1}[1]{\renewcommand{\theoremname}{#1}
    \begin{namedprop}}
    {\end{namedprop}}

\def\R{\mathbb{R}}
\def\E{\mathbb{E}}
\def\P{\mathbb{P}}
\def\N{\mathbb{N}}
\def\eps{\varepsilon}

\begin{document}

\title{Almost Euclidean sections of the N-dimensional cross-polytope
using O(N) random bits}
\author{Shachar Lovett\textsuperscript{1} and Sasha Sodin\textsuperscript{2}}
\footnotetext[1]{Faculty of Mathematics and Computer Science,
The Weizmann Institute of Science, POB 26, Rehovot 76100, Israel.
Email: Shachar.Lovett@weizmann.ac.il. Research supported by
ISF grant 1300/05.}
\footnotetext[2]{School of Mathematics, Raymond and Beverly Sackler
Faculty of Exact Sciences, Tel Aviv University, Tel Aviv, 69978, Israel.
Email: sodinale@tau.ac.il}
\maketitle

\abstract{ It is well known that $\R^N$ has subspaces of dimension
proportional to $N$ on which the $\ell_1$ norm is equivalent to the
$\ell_2$ norm; however, no explicit constructions are known.
Extending earlier work by Artstein--Avidan and Milman, we prove that
such a subspace can be generated using $O(N)$ random bits.}

\section{Introduction}

We study embeddings of $\ell_2$ spaces into $\ell_1$ spaces. Recall
that the $\ell_p$ norm on $\R^N$ is defined by:
\[ \|x\|_p = \left( \displaystyle\sum_{i=1}^N |x_i|^p \right)^{1/p} (p \ge 1)\]

The following inequality holds on $\R^N$:
\[\|x\|_2 \le \|x\|_1 \le \sqrt{N}\|x\|_2 \]

It is well known since the work of Figiel, Lindenstrauss and Milman
\cite{FLM} and Kashin \cite{Kash} that there exists a subspace $E$
of $\R^N$ of dimension $\Theta(N)$ such that for all $x \in E$,
$\|x\|_1 = \Theta(\sqrt{N} \|x\|_2)$ (for the convenience of the
reader, we recall the $\Theta$-notation at the end of the
introduction).

More formally put, for every $0 < \eta < 1$ and every $N \in \N$
(large enough), there exists an $\eta N$-dimensional subspace $E
\hookrightarrow \R^N$ such that for every $x \in E$:
\begin{equation}\label{emb}
c_\eta \sqrt{N} \|x\|_2 \leq \|x\|_1 \leq \sqrt{N} \|x\|_2
\end{equation}
where $c_\eta > 0$ depends only on $\eta$.

The subspace $E$ gives in particular an embedding of $(\R^{\eta N},
\| \cdot \|_2)$ into $(\R^N, \| \cdot \|_1)$. This allows to reduce
various problems in $\ell_2$ norm to corresponding problem in
$\ell_1$ norm, with only a constant blowup in the dimension.

An explicit construction of $E$ would therefore have various
algorithmic applications. This was put forward by Indyk \cite{I,I2},
who proved several related results and applied them to problems in
Computer Science.

No explicit subspace $E$ satisfying (\ref{emb}) has been found so
far (for large $N$). However, it is known that a randomly chosen
subspace, under various natural definitions of distributions of
subspaces, satisfies (\ref{emb}) with probability very close to 1.

In a sense, this situation is typical for various problems in asymptotic
convex geometry, as for numerous properties satisfied by ``random''
high-dimensional objects it is hard to generate a deterministic object
satisfying the property.

To resolve this dissonance, a new line of research was introduced by
Sh.~Artstein-Avidan and V.~Milman. In the innovating work \cite{AM},
the authors proposed to reduce the randomness needed to generate the
random objects. More precisely, they showed that the random constructions
in the proofs of a broad range of theorems, from Milman's Quotient of
Subspace theorem to Zig-Zag approximation, can be performed on the finite
probability space $\{-1, +1\}^R$ equipped with the uniform probability
measure, where $R \in \N$ is reasonably small (the reader may refer
to the work \cite{AM2} by Artstein--Avidan and Milman for further
developments and to the ICM lecture by Szarek \cite{Sz} for a discussion
of these and related issues).

In this case, we say informally that $R$ random bits are used in
the construction. For example, regarding the property (\ref{emb}),
Artstein-Avidan and Milman showed that $O(N \log N)$ random bits
suffice to construct the subspace $E$.

Their proof uses $\eps$-net arguments, and decreasing the number
of random bits beyond $\Omega(N)$ will probably require entirely
new proof ideas. However, the $\log N$ factor in \cite{AM} seemed
to be an artefact of the proof.

In this work, we show that this is indeed the case, and reduce the number
of random bits to $O(N)$ using a modification of the construction
from \cite{AM}.

\begin{thm}\label{thm} For every $0 < \eta < 1$, an $\eta N$-dimensional
subspace of $\R^N$ satisfying (\ref{emb}) can be generated using $O(N)$
random bits. Moreover, the memory needed to generate the subspace
is $O(\log^2 N)$.
\end{thm}

As promised, we recall now the $\Theta$-notation:
\begin{not*}
Let $f,g$ be two functions from $(a, +\infty)$ or $(a, +\infty) \cap
\N$ to $\R_+$. We will write:
\begin{enumerate}
\item $f = O(g)$ if there exist two constants $C > 0$ and $x_0 \geq a$
such that $f(x) \leq C g(x)$ for every $x \geq x_0$;
\item $f = o(g)$ if $f(x)/g(x) \to 0$ as $x \to \infty$;
\item $f = \Omega(g)$ if $g = O(f)$;
\item $f = \omega(g)$ if $g = o(f)$;
\item and finally, $f = \Theta(g)$ if $f = O(g)$ and $f = \Omega(g)$.
\end{enumerate}
\end{not*}

{\em Acknowledgement.} We thank our supervisors, Omer Reingold and
Vitali Milman, for constant support and for their interest in this
work. We are also grateful to Shiri Artstein--Avidan for numerous
discussions and explanations, and in particular for focusing our
attention on bounding the operator norm as the main technical
challenge.


\section{Construction}\label{constr} 

Denote $\xi = 1 - \eta$, $n = \xi N$. We will construct a random $n \times N$
sign matrix $A$ (that is, $A_{ij} = \pm 1$) using $O(N)$ random bits, and
then prove that the kernel
\[ E = \text{Ker} A
    = \left\{ x \in \R^N \, \big| \, Ax = 0 \right\} \]
satisfies (\ref{emb}) with high probability.

Recall the following simple definition:
\begin{dfn}
The Hadamard (or entrywise) product of two $n \times N$ matrices $A_1$
and $A_2$ is the $n \times N$ matrix $A = A_1 \bullet A_2$, defined
by $(A)_{i,j}=(A_1)_{i,j} (A_2)_{i,j}$.
\end{dfn}

Our random matrix $A$ will be the Hadamard product $A_1 \bullet A_2$
of two random matrices $A_1$ and $A_2$, independent of each other.
The construction of $A_1$ and $A_2$ will use two different techniques,
both of them quite common.

\begin{dfn}
A sequence of random variables $X_1,...,X_M$ is called $k$-wise
independent if every $k$ of them are independent.
\end{dfn}

It is well-known that it is possible to construct $M$ $k$-wise
independent random signs from $O(k\log{M})$ truly independent
random signs. More formally, we have:
\begin{lemma1}{Lemma A}
For every $k\leq M$, there exists a subset
\[ \Upsilon_{k,M} \subset \{ -1, 1 \}^M \]
such that $|\Upsilon_{k,M}| = 2^{C_{k,M}}$, $C_{k,M} = O(k \log M)$, and for
the randomly chosen vector $X = (X_1,...,X_M)$ from $\Upsilon_{k,M}$, the
following properties hold:
\begin{enumerate}
\item For $1 \leq m \leq M$,
$\P \{ X_m = -1 \} = \P \{ X_m = 1\} = 1/2$.
\item The coordinates of $X$ are $k$-wise independent.
\item The set $\Upsilon_{k,M}$ is {\em explicit}, meaning that there exists
a bijection $\upsilon_{k,M}: \{-1, 1\}^{C_{k,M}} \to \Upsilon_{k,M}$
that can be computed in time polynomial in $k$ and $M$.
\end{enumerate}
\end{lemma1}

\begin{dfn}\label{def_kind}
The random variables $(X_1, \cdots, X_M)$ satisfying the conditions 1.-2.\ of
Lemma~A are called $k$-wise independent random signs.
\end{dfn}

For completeness, we reproduce a proof of Lemma~A due to Alon, Babai and
Itai \cite{ABI} in Appendix \ref{appind}.

The elements of our first matrix $A_1$ will be $k$-wise independent with
$k = \Theta(\log N)$. That is, $A_1$, regarded as a vector in $\{-1,1\}^{nN}$,
will be a uniformly chosen element of $\Upsilon_{k,nN}$.

\begin{rmk*} Regardless of the distribution of the random sign matrix $A_2$,
the entries $A_{ij}$ of the Hadamard product $A = A_1 \bullet A_2$ are
$k$-wise independent random signs (in the sense of Definition~\ref{def_kind}).
\end{rmk*}

Recall the definition of $\ell_2$ operator norm:
\begin{dfn}
For a matrix A, we define its operator norm as
\[ \|A\| = \max_{x \ne 0} {\frac{\|Ax\|_2}{\|x\|_2}}~.\]
\end{dfn}

The $k$-wise independence of the elements of $A_1$ allows to control the
operator norm of $A$. The following technical lemma may be of independent
interest:
\begin{lemma}\label{opnormlemma}
Let $V$ be any $n \times N$ matrix of $2k$-wise independent random
signs, $k \leq c_2 \sqrt{N}$ (where $c_2 > 0$ is a numerical constant).
Denote $\xi = n/N \le 1$. Then, for $t \geq 0$,
\[\begin{split}
\P{\left\{\frac{1}{\sqrt{N}}\|V\|
    \ge 1 + \sqrt{\xi} + t\right\}}
        &\le 2n \left( 1 + \frac{t}{1+\sqrt{\xi}} \right)^{-2k} \\
        &\le 2n \exp \left\{\frac{-2kt}{1 + \sqrt{\xi} + t}\right\}~.
\end{split}\]
\end{lemma}
We prove the lemma in Section~\ref{opnorm}.

\begin{cor}\label{opnormcor}
Let $0 < \xi < 1$, $n = \xi N$; let $A_1$ be constructed as above
with $k$-wise independent entries, and let $A = A_1 \bullet A_2$,
where $A_2$ is an arbitrary random sign matrix independent of $A_1$.
There exists a numerical constant $C_1>0$ such that for
$k \geq C_1 \log n$,
\[ \P[\|A\|>3\sqrt{N}]<1/n~.\]
\end{cor}

We now head to construct a probability space for $A_2$; we use random
walks on expander graphs (see Hoory, Linial and Wigderson \cite{HLW}
for an extensive survey). Let us recall the basic definitions.

Let $G = (\mathcal{V}, \mathcal{E})$ be a $d$-regular graph; the value of $d$ plays no
significant role in the estimates, so the reader may assume $d = 4$.
Let $P^G$ be the transition matrix of the random walk of $G$:
\[ P^G_{uv} = \begin{cases}
1/d, &(u,v) \in \mathcal{E} \\
0, &(u,v) \notin \mathcal{E}.
\end{cases}\]
Denote by $1 = \lambda_1 \geq \lambda_2 \geq \lambda_3 \geq \cdots$
the eigenvalues of $P^G$ arranged in decreasing order, and denote
$\lambda = \max_{i \geq 2} |\lambda_i|$.

In this notation, the graph $G$ is called a $(|\mathcal{V}|, d, \lambda)$-graph.
We will only need the following fact (cf.\ \cite{HLW}, \cite{AM}):

\begin{fact*} For any $d\geq 3$ and any number of vertices $|\mathcal{V}|$
(big enough), there exists a $(|\mathcal{V}|, d, \lambda)$-graph
$G = (\mathcal{V} = \{1, 2, \cdots, |\mathcal{V}|\}, \mathcal{E})$ such that
\begin{enumerate}
\item $\lambda < 0.95$ and
\item $G$ is {\em explicit}, formally meaning that set of neighbours
\[ \{u \in \mathcal{V} \, | \, (u, v) \in \mathcal{E}\} \]
of any vertex $v \in \mathcal{V}$ can be computed in time that is polynomial
in $\log |\mathcal{V}|$.
\end{enumerate}
\end{fact*}

Sometimes we will call such a graph an expander graph with parameter $\lambda$.

Let $G = (\mathcal{V}, \mathcal{E})$ be an expander graph, with vertices $\mathcal{V}$
indexed by the elements of $\Upsilon_{4,N}$. Let $v_1, v_2, \cdots, v_n$ be a
random walk of length $n$ in $G$, starting from a random element of $\mathcal{V}$.
Write the sign vectors corresponding to $v_1,\cdots,v_n$ in $\Upsilon_{4,N}$ as
the rows of $A_2$.

The use of expander graphs is similar to \cite{AM}; however, we use
constant degree expanders. We also show it suffices to use 4-wise
independent rows rather than truly independent rows. This enables
the computation to be performed using less memory ($O(\log^2{N})$).

Note that the construction uses in total
\begin{equation}\label{nbits}
\begin{split}
O(\log n \log (Nn)) &+ O(\log N) + O(n \log d) \\
    &= O(n + \log n \log N) = O(N)
\end{split}
\end{equation}
random bits.  Also, we have the following:

\begin{lemma}\label{singveclemma}
Let $A_1$ be any constant sign matrix, and let $A_2$ be constructed as above.
For every $x \in \R^N$ and any $\eps \leq c_\lambda \sqrt{\xi}$,
\[\P \left\{\|Ax\|_2 < 6 \eps \sqrt{N} \|x\|_2\right\}
    < C_\lambda p_\lambda^n~,\]
where the constants $C_\lambda, c_\lambda > 0$ and $0 < p_\lambda < 1$ depend on the
parameter $\lambda \in [0, 1)$ of the graph $G$.
\end{lemma}

\begin{cor}\label{singleveccor}
The statement of the lemma remains true if we change $A_1$ from constant to drawn
from any distribution.
\end{cor}

We prove this lemma in Section~\ref{singvec}; the proof is a variation on the ideas
from Artstein-Avidan and Milman \cite{AM}.

Now we can reformulate our main result.

\begin{thm}\label{thm2} Let $A_1$ and $A_2$ be constructed as above ($A_1$ has
$\Theta(\log n)$ independent entries, the rows of $A_2$ come from a random
walk on an expander); let $A = A_1 \bullet A_2$, $E = \text{Ker} \, A$.
Then, with probability $1 - o(1)$,
\begin{equation}\label{emb1}
\frac{c' \xi}{\sqrt{\log 1/\xi}} \sqrt{N} \|x\|_2
    \leq \|x\|_1 \leq \sqrt{N} \|x\|_2
        \quad \text{for every} \quad x \in E~,
\end{equation}
where $c'>0$ is a universal constant.
\end{thm}

The proof uses the Lemmata formulated above as well as the following standard
lemma from asymptotic convex geometry.

\begin{lemma1}{Lemma B}
Let A be a random $n \times N$ sign matrix such that:
\begin{enumerate}
\item $\P[\|A\| > 3\sqrt{N}] \le q$;
\item There exist $0<p<1$, $\eps>0$ and $C > 0$ such that for every $y \in \R^N$,
\[ \P \left\{ \|Ay\|_2 < 6\eps \sqrt{N} \|y\|_2 \right\}< C p^n~.\]
\end{enumerate}

Then with probability at least
\[ 1-q-p^{\Theta(n)}\]
over the choice of $A$, we have:
\[ \|x\|_1 \ge \delta \sqrt{N} \|x\|_2
    \quad \text{for every} \quad
    x \in \text{Ker} A~,\]
where we can take
\[ \delta = \frac{c\eps}{\sqrt{\frac{1}{\xi}\log{\frac{1}{p}}
    \log{(\frac{1}{\xi} \log{\frac{1}{p}})}}} ~,\]
$c>0$ being a universal constant.
\end{lemma1}

For completeness, we prove Lemma~B in Appendix~\ref{agsect}.

\begin{proof}[Proof of Theorem~\ref{thm2}]
According to Corollary~\ref{opnormcor} the random matrix $A$ satisfies the condition 1.\
of Lemma~B with $q = 1/n$. According to Corollary~\ref{singleveccor} $A$ also satisfies 2.,
with $p = p_\lambda$, $C = C_\lambda$ and $\eps = c_\lambda \sqrt{\xi}$. Now apply Lemma~B;
note that $\lambda \leq 0.95 < 1$ is bounded away from $1$ and hence $p_\lambda$ and
$C_\lambda$ may be replaced by universal constants ($p_{0.95}$ and $C_{0.95}$, resp.)
\end{proof}

Clearly, Theorem~\ref{thm2} implies Theorem~\ref{thm}.

\section{Operator norm of a matrix with $2k$-wise independent entries}\label{opnorm} 

\begin{proof}[Proof of Lemma~\ref{opnormlemma}]
We start by bounding the expectation of $\|V\|^{2k}$. For a real symmetric
$n \times n$ matrix $W$, denote by $\lambda_1(W), \cdots, \lambda_n(W)$
the eigenvalues of $W$, and let $\lambda_{\max}(W) = \max_i \lambda_i(W)$.
Observe that
\[ \|V\|^2 = \lambda_{\max}(V^t V) = \lambda_{\max}(V V^t) \]
and hence:

\[\begin{split}
\E \|V/\sqrt{N}\|^{2k} &= \E \lambda_{\max}(V V^t/N)^k \\
&\le \E \displaystyle\sum_{i=1}^n \lambda_i(V V^t/N)^k
= \E \text{Tr}((V V^t/N)^k)~.
\end{split}\]

The trace of $(V V^t)^k$ is equal to
\[ \sum V_{i_1,j_1}V_{i_2,j_1}V_{i_2,j_2} V_{i_3,j_2}
    \cdots V_{i_k,j_k}V_{i_1,j_k}~,\]
where the sum is over closed paths $(i_1,j_1,...,i_k,j_k,i_1)$ in the
bipartite graph $K_{n,N}$. The expectation of each term in the sum is $0$
if there is some $V_{i,j}$ that appears an odd number of times, and $1$
if all the terms appear an even number of times. So, the expectation
is equal to the number $m(k;n,N)$ of closed even paths of length $2k$
in $K_{n,N}$, starting on the side of size $n$ (an even path is a path
in which every edge appears an even number of times).

Instead of estimating this expectation directly, we follow an idea
of Aubrun \cite{Au} and take a different route. The trace of $(V V^t)^k$
is a sum over products of powers of at most $2k$ elements from V, and so,
since the elements of $V$ come from a $2k$-wise independent probability space,
the expectation is the same as if the elements of $V$ were truly
independent. Hence, we can use estimates known for matrices with
i.i.d.\ elements.

We chose to use such an estimate for matrices with Gaussian i.i.d
elements. Let $\widetilde{V}$ be an $n \times N$ matrix, whose entries
are independent, $\widetilde{V}_{i,j} \sim N(0,1)$.
For every entry $1 \leq i \leq n$, $1 \leq j \leq N$ and every
integer $l \geq 1$ we have:
\[ \E \widetilde{V}_{i,j}^{2l} \ge (\E \widetilde{V}_{i,j}^2)^{l}
    = 1 = \E V_{i,j}^{2l}\, ; \quad
\E \widetilde{V}_{i,j}^{2l+1} = 0 = \E V_{i,j}^{2l+1}~.\]

Therefore
\[\begin{split}
\E \text{Tr}((V V^t/N)^k)
    &\le \E \text{Tr}((\widetilde{V} \widetilde{V}^t/N)^k)
    = \E \displaystyle\sum_{i=1}^n \lambda_i(\widetilde{V} \widetilde{V}^t/N)^k \\
    &\le n \E \lambda_{\max}(\widetilde{V} \widetilde{V}^t/N)^k
    = n \, \E \| \widetilde{V}/\sqrt{N}\| ^{2k}~.
\end{split}\]

We use the following bound for Gaussian random matrices with independent
entries (see Davidson--Szarek \cite[Thm.~II.13]{DS}, extending an idea of
Y.~Gordon):
\[ \P\left\{\|\widetilde{V}/\sqrt{N}\| \ge 1 + \sqrt{\xi} + t\right\} <
    \exp(-Nt^2/2)~, \quad t \geq 0~.\]

Now,
\begin{multline*}
\E \|\widetilde{V}/\sqrt{N}\|^{2k}
    = \displaystyle\int_{0}^{\infty} 2k t^{2k-1}
        \P \left\{\|\widetilde{V}/\sqrt{N}\| \ge t \right\}dt \\
    < (1+\sqrt{\xi})^{2k}
        + 2k \displaystyle\int_0^{\infty} (1+\sqrt{\xi}+u)^{2k-1}
            \exp(-Nu^2/2)\, du~.
\end{multline*}
It is easy to see that the second term is smaller than the first one:
\begin{multline*}
2k \int_0^{\infty} (1+\sqrt{\xi}+u)^{2k-1} \exp(-Nu^2/2)\, du \\
< 2k (1+\sqrt{\xi})^{2k-1} \int_0^\infty
    \exp \left\{ \frac{2k-1}{1+\sqrt{\xi}} \, u - Nu^2/2 \right\} du \\
< \frac{2k}{\sqrt{N}} (1 + \sqrt{\xi})^{2k-1} \int_{-\infty}^\infty
    \exp \left\{ \frac{2k-1}{\sqrt{N} + \sqrt{n}} \, u - u^2/2 \right\} du \\
= (1 + \sqrt{\xi})^{2k-1} \frac{\sqrt{8\pi} \, k}{\sqrt{N}}
    \exp \left\{ \frac{1}{2} \left( \frac{2k-1}{\sqrt{N}+\sqrt{n}} \right)^2 \right\} \\
= (1 + \sqrt{\xi})^{2k} \times O(k/\sqrt{N}) \times e^{O(k^2/N)}~.
\end{multline*}
If $k \leq c_2 \sqrt{N}$ (for an appropriately chosen numerical constant $c_2>0$),
the product of the $O$-terms is not greater than $1$. Hence
\[ \E \|\widetilde{V}/\sqrt{N}\|^{2k} < 2(1+\sqrt{\xi})^{2k}~, \]
implying that
\[ \E \|V/\sqrt{N}\|^{2k} < 2n(1+\sqrt{\xi})^{2k}~.\]
Now by Chebyshev's inequality
\[\P\left\{ \|V/\sqrt{N}\| \ge 1+ \sqrt{\xi} + t\right\}
    \le \frac{\E\|V/\sqrt{N}\|^{2k}}{(1+\sqrt{\xi}+t)^{2k}}
    < 2n \left( \frac{1+\sqrt{\xi}}{1+\sqrt{\xi}+t} \right)^{2k} \]
\end{proof}

\begin{rmks*}\hfill
\begin{enumerate}
\item The lemma shows that for $k = \Omega(\log N)$ the operator
norm of $V/\sqrt{N}$ is not much larger than $1 + \sqrt{\xi}$. This
matches the bound for matrices with independent entries
(cf.\ Geman \cite{G}).
\item A more direct proof would be to bound the numbers $m(k;n,N)$
directly, as in the work of Geman \cite{G}. This would yield
an estimate similar to the one we get.
\end{enumerate}
\end{rmks*}

\section{Bound for a single vector}\label{singvec} 

Fix $x$, $\| x \|_2 = 1$; let us bound the probability
\[ \P \left\{ \| A x \|_2 < 6 \eps \sqrt{N} \right\}\]
when $A = A_1 \bullet A_2$, $A_1$ is a fixed sign matrix and
$A_2$ is generated from a random walk on an expander as
explained in Section~\ref{constr}.

Recall that $G = (\mathcal{V}, \mathcal{E})$ is a $d$-regular
graph with $2^{O(\log N)}$ vertices, and $P^G$ is the transition
matrix of the random walk on $G$; $\lambda$ is the second largest
absolute value of an eigenvalue of $P^G$.

First we bound from below the probability that a
coordinate of $Ax$ is not very small.

\begin{lemma}\label{l1}
Let $\Psi$ be a random vector in $\{ -1, +1\}^N$ with 4-wise independent
coordinates. Then
\[ \P \left\{ \langle \Psi, x \rangle^2 \geq 1/2 \right\} \geq 1/12~. \]
\end{lemma}

\begin{proof}
First,
\[ \E \langle \Psi, \, x \rangle^2
    = \sum_{i,j=1}^N x_i x_j \E \Psi_i \Psi_j
    = \sum_{i=1}^N x_i^2 = 1~; \]
\[\begin{split}
\E \langle \Psi, \, x \rangle^4
    &= \sum_{i,j,k,l=1}^N x_i x_j x_k x_l \, \E \Psi_i \Psi_j \Psi_k \Psi_l \\
    &= \sum_{i=1}^N x_i^4 + 6 \sum_{1 \leq i < j \leq N} x_i^2 x_j^2
    < 3 \left( \sum_{i=1}^N x_i^2 \right)^2 = 3~.
\end{split}\]

Recall the Paley--Zygmund inequality \cite{PZ}:

\begin{lemma*}[Paley--Zygmund] If $Z \geq 0$ is a random variable with
finite second moment, $0 < \theta < 1$, then
\[ \P \left\{ Z \geq \theta \E Z \right\}
    \geq (1-\theta)^2 \frac{\E(Z)^2}{\E(Z^2)}~. \]
\end{lemma*}

Applying the inequality for $Z = \langle \Psi, \, x \rangle^2$, $\theta = 1/2$,
we obtain the statement of the lemma.

\end{proof}

\begin{proof}[Proof of Lemma~\ref{singveclemma}]
Let us show that a constant fraction of the rows $\psi_i$ of $A$
satisfy w.h.p
\begin{equation}\label{sr}
\langle \psi_i, x \rangle \geq 1/2~.
\end{equation}
For fixed $A_1$ and $1 \leq i \leq n$, the coordinates of $\psi_i$ are $4$-wise
independent; therefore by Lemma~\ref{l1} there is a subset $S_i \subset \mathcal{V}$
such that $|S_i|/|\mathcal{V}| \geq 1/12$, and the $i$-th $\psi_i$ of $A$ satisfies
(\ref{sr}) iff the $i$-th row $v_i$ of $A_2$ lies in $S_i$.

We need a modification of Kahale's Chernoff-type bound on expanders
\cite{K}, see also Alon, Feige, Wigderson and Zuckerman
\cite[Theorem 4]{AFWZ}, Art\-stein-Avi\-dan and Milman \cite[Section
4]{AM}, and Hoory, Linial and Wigderson \cite[Theorem 3.11]{HLW} for
related results\footnote{Added in proof: an even stronger result was
recently proved. See theorem 5.4 in E. Mossel, R. O'Donnell, O.
Regev, J. Steif and B. Sudakov, Non-Interactive Correlation
Distillation, Inhomogeneous Markov Chains and the Reverse
Bonami-Beckner Inequality, Israel Journal of Mathematics 154 (2006),
299-336.}

\begin{lemma}\label{chl} Let $G = (\mathcal{V}, \mathcal{E})$ be a graph;
as before, let $1 = \lambda_1 \geq \lambda_2 \geq \lambda_3 \geq \cdots$
be the eigenvalues of $P^G$; denote $\lambda = \max_{i \geq 2} |\lambda_i|$.
The probability that a random walk on $G$, starting from a random point in $\mathcal{V}$,
is in $S_i$ on the $i$-th step, $i = 1, 2, \cdots, k$, is at most
\[ \prod_{i=1}^{k-1}
    \sqrt{\lambda + (1-\lambda) \frac{|S_i|}{|V|}} \,
    \sqrt{\lambda + (1-\lambda) \frac{|S_{i+1}|}{|V|}}~.\]
\end{lemma}

\begin{proof}[Proof of Lemma~\ref{chl}]
Denote $e = (1, 1, \cdots, 1)/\sqrt{|\mathcal{V}|}$, and denote
by $\Pi_i$ the projector on the coordinates in $S_i$. Then the probability
in question equals
\begin{multline}
\langle \Pi_k P^G \Pi_{k-1} P^G \cdots P^G \Pi_1 e, e \rangle \\
    \leq \| \Pi_k P^G \Pi_{k-1} \| \times \| \Pi_{k-1} P^G \Pi_{k-2} \| \times
        \cdots \times \| \Pi_2 P^G \Pi_1 \|~,
\end{multline}
where we used the submultiplicativity of operator norm and the equality
$\Pi_i^2 = \Pi_i$. Let us bound the norms
\[ \| \Pi_{i+1} P^G \Pi_i \| = \max_{\|g\|_2 = 1}  \| \Pi_{i+1} P^G \Pi_i g\|_2~.\]

First of all, the vector $g$ for which the maximum is attained is supported
in $S_i$; hence $\Pi_i g = g$. Let us decompose $g = \alpha e + \beta v$,
where $\alpha^2 + \beta^2 = 1$ and $v$ is a unit vector orthogonal to $e$.

Note that
\[ |\alpha| = |\langle g, e \rangle|
    \leq \|g \|_1/\sqrt{|\mathcal{V}|}
    \leq \sqrt{\frac{|S_{i}|}{|\mathcal{V}|}} \,  \|g\|_2
    = \sqrt{\frac{|S_{i}|}{|\mathcal{V}|}}~. \]
Therefore $P^Gg = \alpha e + \beta P^Gv$. Now,
\[ \| \Pi_{i+1} P^Gg \|_2
    = \max_{\|h\|_2 = 1} \langle \Pi_{i+1} P^Gg, h \rangle
    = \max_{\|h\|_2 = 1} \langle P^Gg, \Pi_{i+1} h \rangle~;\]
we may assume that $h$ is supported in $S_{i+1}$. Let $h = \alpha' e + \beta'v'$,
where $v'$ is a unit vector orthogonal to $e$; as before,
\[ \alpha'^2 + \beta'^2 = 1 \quad \text{and} \quad
    |\alpha'| \leq \sqrt{\frac{|S_{i+1}|}{|\mathcal{V}|}}~.\]
Hence
\[\begin{split}
\langle P^Gg, h \rangle
    &= \alpha \alpha' + \beta \beta' \langle P^Gv, v' \rangle
    \leq \alpha \alpha' + \lambda \beta \beta' \\
    &\leq \sqrt{\alpha^2 + \lambda \beta^2} \,
        \sqrt{\alpha'^2 + \lambda \beta'^2} \\
    &= \sqrt{\lambda + (1-\lambda) \alpha^2} \,
        \sqrt{\lambda + (1-\lambda) \alpha'^2} \\
    &\leq \sqrt{\lambda + (1-\lambda) \frac{|S_i|}{|\mathcal{V}|}} \,
        \sqrt{\lambda + (1-\lambda) \frac{|S_{i+1}|}{|\mathcal{V}|}}~.
\end{split}\]
\end{proof}

Now, if $\|A x\|_2 < 6 \eps \sqrt{N}$, $A$ has at most $72 \eps^2 N$ rows $\psi$
such that
\[ \langle \psi, x \rangle^2 \geq 1/2~.\]
By Lemma~\ref{chl}, the probability of this event is at most
\begin{multline}
\binom{n}{[72 \eps^2 N]}
        \left( \frac{11}{12}(1 - \lambda) + \lambda \right)^{n - [72 \eps^2 N] - 1} \\
    \leq 2 \left(\frac{e\xi}{72 \eps^2}\right)^{72n\eps^2/\xi}
        \left( \frac{11}{12}(1 - \lambda) + \lambda \right)^{n - 72 n \eps^2/ \xi}~.
\end{multline}

For $\eps$ small enough, this probability is exponentially small. More formally, it
is easy to see that there exist some constants $C_\lambda \geq 1 > c_\lambda >0$
and $0 < p_\lambda < 1$ depending only on $\lambda$, such that
\begin{equation}\label{singlevec}
\P \left\{ \|A x\|_2 < 6 \eps \sqrt{N} \right\} \leq C_\lambda p_\lambda^n
    \quad \text{if} \quad 0 < \eps \leq c_\lambda \sqrt{\xi}~.
\end{equation}

Lemma~\ref{singveclemma} is proved.

\end{proof}

\appendix

\section{Construction of $k$-wise independent random bits}\label{appind} 

For completeness, we recall the construction of $2^r-1$ $k$-wise
independent random bits from $kr$ independent random bits due to Alon, Babai
and Itai \cite{ABI}. It will be more convenient to work with vectors
of $\{0,1\}$ rather than $\{-1,+1\}$.

Let
\[ \alpha_1, \cdots, \alpha_{2^r-1} \in \text{GF}(2^r) \]
be the non-zero elements of the finite field of cardinality $2^r$.
$\text{GF}(2^r)$ is a linear space over $\text{GF}(2)$; hence we may represent
an element $\alpha \in \text{GF}(2^r)$ as an $r$-tuple
$\widetilde{\alpha} \in \text{GF}(2)^r$.

Consider the matrix
\[ M = \left(\begin{matrix}
    1 & \alpha_1 & \alpha_1^2 &\cdots &\alpha_1^{k-1}\\
    1 & \alpha_2 & \alpha_2^2 &\cdots &\alpha_2^{k-1}\\
    \cdots & \cdots &\cdots & \cdots & \cdots\\
    1 & \alpha_{2^r-1} & \alpha_2^2 &\cdots &\alpha_{2^r-1}^{k-1}
\end{matrix}\right)~.\]

Every $k$ rows of $M$ form a Van der Monde matrix, and in particular
are linearly independent. Let

\[ \widetilde{M} = \left(\begin{matrix}
    1 & \widetilde{\alpha}_1 & \widetilde{\alpha}_1^2
        &\cdots &\widetilde{\alpha}_1^{k-1}\\
    1 & \widetilde{\alpha}_2 & \widetilde{\alpha}_2^2
        &\cdots &\widetilde{\alpha}_2^{k-1}\\
    \cdots & \cdots &\cdots & \cdots & \cdots\\
    1 & \widetilde{\alpha}_{2^r-1} & \widetilde{\alpha}_2^2
        &\cdots &\widetilde{\alpha}_{2^r-1}^{k-1}
\end{matrix}\right)\]
be the corresponding $kr \times (2^r-1)$ matrix over $\text{GF}(2)$;
its rows are also linearly independent. Now let $Z$ be a random
vector distributed uniformly in $\text{GF}(2)^{kr}$; let
$X = \widetilde{M} Z$.

\begin{claim*} The coordinates of the vector $X$ are $k$-wise independent.
\end{claim*}

\begin{proof}
For every set of indices
$\emptyset \neq I \subset \{1, \cdots, 2^r-1\}$
such that $|I| = k$, the matrix $\widetilde{M}_I$ formed from the
corresponding rows of $\widetilde{M}$ is of rank $k$; that is,
$\widetilde{M}_I$ is surjective and the preimages of the vectors
in $\{0,1\}^k$ are of equal size. The vector $Z$ is distributed
uniformly in $\text{GF}(2)^{kr}$; hence the vector
$\left(X_i\right)_{i \in I} = \widetilde{M}_I Z$ is uniformly
distributed in $\text{GF}(2)^k$.
\end{proof}

\section{Proof of Lemma~B}\label{agsect}

The proof of Lemma~B is based on $\eps$-net arguments.

\begin{dfn}
Let $S \subset \R^N$ be a convex set. A (finite) subset
$\mathcal{N} \subset S$ is called an $\eps$-net in $S$
if for every $x \in S$ there exists $y \in \mathcal{N}$
such that $\|x-y\|_2 \leq \eps$.
\end{dfn}

\begin{not*} Let $t > 0$ and let $K \subset \R^n$ be a convex body.
As usual, denote
\[ tK = \{tx \, | \, x \in K \}~.\]
\end{not*}

Similarly to \cite{AM}, we use the following result,
due to Sch\"utt \cite{Sch}:
\begin{thm*}[Sch\"utt]
The exists a universal constant $c>0$ such that for any $\zeta > 0$
and $\theta \geq c \sqrt{\frac{1}{\zeta} \log \frac{1}{\zeta}}$ there
exists a $\theta$-net $\mathcal{N}$ in $\sqrt{N} B_1^N$ such that
$|\mathcal{N}| \leq e^{\zeta N}$.
\end{thm*}

\begin{proof}[Proof of Lemma B]
Pick $0 < \zeta < \xi \log \frac{1}{p}$; then $e^\zeta < 1/p^\xi$.
Set
\[ \delta = \frac{\eps}{c \sqrt{\frac{1}{\zeta} \log \frac{1}{\zeta}}}~.\]
Scaling the result of Sch\"utt's theorem times $\delta$, we get
an $\eps$-net $\mathcal{N}$ in $\delta \sqrt{N} B_1^N$,
$|\mathcal{N}| \leq e^{\zeta N}$.

By our assumptions, for every $y \in \mathcal{N}$
\[ \P \left\{ \|Ay\|_2 < 6\eps \sqrt{N} \|y\|_2 \right\} < C p^n~,\]
and so the probability that there exists $y \in \mathcal{N}$ with
\[ \|Ay\|_2 < 6\eps \sqrt{N} \|y\|_2 \]
is at most
\[ C e^{\zeta N} p^n = p^{\Theta(n)}~.\]

Assume that for every $y \in \mathcal{N}$ we have
\[ \|Ay\|_2 \ge 6\eps \sqrt{N} \|y\|_2~,\]
and also that $\| A \| \le 3\sqrt{N}$. This event happens with
probability at least $1-q-p^{\Theta(n)}$. We will show that whenever
these two conditions hold, every $x \in \text{Ker} A$ satisfies
\[ \|x\|_1 \ge \delta \sqrt{N} \|x\|_2~. \]
It is enough to show this for $x$ with $\|x\|_2=1$.

Take any $x \in \R^N$ with $\|x\|_1 < \delta \sqrt{N}$ and
$\|x\|_2=1$. We will show $x \notin \text{Ker}(A)$. First,
$x \in \delta \sqrt{N} B_1^N$, and so there exists $y \in \mathcal{N}$
such that $\|x - y \| \leq \eps$. Now we have:

\[\begin{split}
\|Ax\|_2 &\ge \|Ay\|_2 - \|A(x-y)\|_2
    \ge 6\eps\sqrt{N}\|y\|_2 - \|A\| \|x-y\|_2 \\
    &\ge 6\eps (1-\eps) \sqrt{N} - 3\eps\sqrt{N} > 0~,
\end{split}\]

where we used the fact that
\[ \|y\|_2 \ge \|x\|_2 - \|x-y\|_2 \ge 1-\eps~.\]
\end{proof}

\end{document}